\documentclass[11pt]{amsart}
\usepackage{footnote}
\newtheorem{theorem}{Theorem}[section]
\newtheorem{corollary}[theorem]{Corollary}
\newtheorem{lemma}[theorem]{Lemma}
\newtheorem{proposition}[theorem]{Proposition}

\theoremstyle{definition}
\newtheorem{definition}[theorem]{Definition}
\theoremstyle{remark}
\newtheorem{remark}[theorem]{Remark}
\numberwithin{equation}{section}

\newcommand{\CC}{\mathcal C}

\newcommand{\C}{\mathbb C}
\newcommand{\R}{\mathbb R}

\begin{document}

\title{Plurisubharmonic exhaustion functions and almost complex Stein structures}

\author{Klas Diederich and Alexandre Sukhov}

\address{\begin{tabular}{lll}
Klas Diederich & & Alexandre Sukhov\\
Department of Mathematics & & U.S.T.L. \\
University of Wuppertal & & Cit\'e Scientifique \\
D-42097 Wuppertal & & 59655 Villeneuve d'Ascq Cedex\\
GERMANY & & FRANCE\\
 & & \\
{\rm Klas.Diederich@math.uni-wuppertal.de} & &
{\rm sukhov@math.univ-lille1.fr}
\end{tabular}
}

\subjclass[2000]{32H02, 53C15 \smallskip\\ \hspace*{3mm} Submitted March 2006}


\begin{abstract}
We prove that a relatively compact pseudoconvex domain with smooth  boundary in
an almost complex manifold admits a bounded strictly plurisubharmonic
exhaustion function. We use this result in order to study
convexity and hyperbolicity
properties  of these domains and the contact geometry of their boundaries.
\end{abstract}

\maketitle

\section{Introduction}
One of the fundamental results of classical complex analysis establishes
the equivalence between the holomorphic disc convexity of a domain in an
affine complex space,    the Levi convexity (the positive semidefiniteness of the Levi form)
of its boundary and the existence of a strictly
plurisubharmonic  exhaustion function. On the other hand, in the works of
 Y.Eliashberg and M.Gromov, D.McDuff \cite{GrEl, Mc} and other authors the
convexity properties of strictly pseudoconvex domains in almost complex
manifolds are substantially used and allow to obtain  many interesting results concerning symplectic and
contact structures.  It turned out that  the notion of pseudoconvexity
playing a fundamental role in classical complex analysis
 admits deep analogs in the symplectic
category. Furthermore,
Y.Eliashberg and W.Thurston \cite{ElTh} introduced and
studied  so-called confoliations which can be viewed as contact structures
with degeneracies. Their theory, developed in \cite{ElTh}
mainly for
manifolds of real dimension 3,  links the geometry and topology of contact
structures with the theory of foliations. One of the main examples of
confoliations is given by  the distribution  of holomorphic tangent spaces on
the boundary of a
  weakly pseudoconvex domain in an almost complex manifold (in general,
 with a non-integrable almost complex structure). However, in contrast
 to the situation for the classical case of  $\C^n$,  only quite a few
 properties of such domains are known.  The goal of our paper is to study
the convexity properties of (weakly) pseudoconvex  domains in almost complex
manifolds.

The paper is organized as follows. Section 2 and 3 are essentially
preliminary and contain the properties of almost complex structures used in
the proofs of main results.  In section 4 we prove that if
 $(M,J)$ is an almost complex manifold admitting a strictly plurisubharmonic
 function and $\Omega$ is a relatively compact domain in $M$ with smooth
 boundary such that the Levi form of $b\Omega$ is positive semidefinite at
 every point, then $\Omega$ admits a  bounded strictly  plurisubharmonic
 exhaustion function (Theorem \ref{theo1}). This generalizes the known
 result of K.Diederich - J.E.Fornaess \cite{DiFo1,DiFo2} dealing with the case of
 domains in $\C^n$. In particular, this means that the domains satisfying the
 hypothesis of Theorem \ref{theo1}   are Stein manifolds in the sense
 of Y.Eliashberg - M.Gromov
 \cite{GrEl} and admit the canonical symplectic structure defined by the Levi
 form of a  strictly plurisubharmonic exhaustion function.\smallskip\\
As an application we obtain a characterization of
 pseudoconvex
domains in almost complex manifolds similar to the classical results of
complex analysis (Theorem \ref{theo2}).\smallskip\\
As another consequence of Theorem \ref{theo1}  we obtain
that in its setting the domain $\Omega$ is taut, that is, every sequence $(D_k)$
of
holomorphic discs in $\Omega$ either contains a subsequence, convergent in the
compact open topology, or $(D_k)$ is compactly divergent (Theorem \ref{theo3}).\smallskip\\
Finally, we prove a result on the approximation of confoliations (in the sense of
Y.Eliashberg - W.Thurston \cite{ElTh}) by contact structures (Theorem \ref{theo4}).

\section{Almost complex manifolds}

All manifolds and almost complex structures are supposed to be of class $C^\infty$
though the main results require only a lower regularity.
Let $(\tilde M, \tilde J)$ and $(M,J)$ be almost complex manifolds and let $f$ be
a smooth map from $\tilde M$ to $M$. We say that $f$ is
 $(\tilde J,J)$-holomorphic if $df \circ \tilde J = J \circ df$.
 Let $\mathbb D$ be the unit
disc in $\C$ and $J_{st}$ be the standard  structure on $\C^n$
for every $n$. If $(M',J')=(\mathbb D,J_{st})$,
we call $f$ a  $J$-holomorphic disc in $M$.

 Every almost complex manifold
$(M,J)$ can be viewed locally as the unit ball $\mathbb B$ in
$\mathbb C^n$ equipped with a small almost complex
deformation of $J_{st}$. Indeed, we have the following frequently used statement.
\begin{lemma}
\label{lemma1}
Let $(M,J)$ be an almost complex manifold. Then for every point $p \in
M$, every real $\alpha \geq 0$ and   $\lambda_0 > 0$ there exist a neighborhood $U$ of $p$ and a
coordinate diffeomorphism $z: U \rightarrow \mathbb B$ such that
$z(p) = 0$, $dz(p) \circ J(p) \circ dz^{-1}(0) = J_{st}$  and the
direct image $ z_*(J): = dz \circ J \circ dz^{-1}$ satisfies $\vert\vert z_*(J) - J_{st}
\vert\vert_{\CC^\alpha(\bar {\mathbb B})} \leq \lambda_0$.
\end{lemma}
\proof There exists a diffeomorphism $z$ from a neighborhood $U'$ of
$p \in M$ onto $\mathbb B$ satisfying $z(p) = 0$ and $dz(p) \circ J(p)
\circ dz^{-1}(0) = J_{st}$. For $\lambda > 0$ consider the dilation
$d_{\lambda}: t \mapsto \lambda^{-1}t$ in $\R^{2n}$ and the composition
$z_{\lambda} = d_{\lambda} \circ z$. Then $\lim_{\lambda \rightarrow
0} \vert\vert (z_{\lambda})_{*}(J) - J_{st} \vert\vert_{\CC^\alpha(\bar
{\mathbb B})} = 0$ for every real $\alpha \geq 0$. Setting $U = z^{-1}_{\lambda}(\mathbb B)$ for
$\lambda > 0$ small enough, we obtain the desired statement. \qed

\bigskip

In the sequel we often denote $z_*(J)$ just by $J$ when local coordinates are fixed.
Let $(M,J)$ be an almost complex manifold. We denote by $TM$ the real
tangent bundle of $M$ and by $T_\C M$ its complexification. Recall
that $T_\C M = T^{(1,0)}M \oplus T^{(0,1)}M$ where
$T^{(1,0)}M:=\{ V \in T_\C M : JV=iV\} = \{X -iJ X, X \in
TM\},$
and $T^{(0,1)}M:=\{ V \in T_\C M : JV=-iV\} = \{X +
iJ X, X \in TM\}$.
 Let $T^*M$ denote the cotangent bundle of  $M$.
Identifying $\C \otimes T^*M$ with
$T_\C^*M:=Hom(T_\C M,\C)$ we define the set of complex
forms of type $(1,0)$ on $M$ by~:
$
T^*_{(1,0)}M=\{w \in T_\C^* M : w(V) = 0, \forall V \in T^{(0,1)}M\}
$
and the set of complex forms of type $(0,1)$ on $M$ by~:
$
T^*_{(0,1)}M=\{w \in T_\C^* M : w(V) = 0, \forall V \in T^{(1,0)}M\}
$.
Then $T_\C^*M=T^*_{(1,0)}M \oplus T^*_{(0,1)}M$.
This allows to define the operators $\partial_J$ and
$\bar{\partial}_J$ on the space of smooth functions defined on
$M$~: given a complex smooth function $u$ on $M$, we set $\partial_J u =
du_{(1,0)} \in T^*_{(1,0)}M$ and $\bar{\partial}_Ju = du_{(0,1)}
\in T^*_{(0,1)}M$. As usual,
differential forms of any bidegree $(p,q)$ on $(M,J)$ are defined
by means of the exterior product.

In what follows we will often work in local coordinates. Fixing local
coordinates
on a connected open neighborhood of a point $p \in (M,J)$, we can view it
 as a neighborhood $U$ of the origin (corresponding to $p$) in $\R^{2n}$ with the
 standard coordinates $(x,y)$. After an additional linear transformation we
 may assume that $J(0) = J_{st}$; recall that
\begin{eqnarray*}
J_{st}(\frac{\partial }{\partial x}) = \frac{\partial }{\partial y},
J_{st}(\frac{\partial }{\partial y}) = -\frac{\partial }{\partial x}.
\end{eqnarray*}
In these coordinates the structure $J$ can be viewed as a smooth real $(2n \times 2n)$-matrix
function $J:U \longrightarrow M_{2n}(\R)$ satisfying $J^2 = -I$ (the $2n$ unit
matrix).
After the above complexification we can identify $\R^{2n}$ and $\C^n$ with the
standard complex coordinates $z = x + iy$ and so deal with $T^{(1,0)}(M)$
instead of $T(M)$.

Consider a $J$-holomorphic disc
$z: \mathbb D \longrightarrow U$. After a straightforward verification  the
conditions $J(z)^2 = -I$ and $J(0) = J_{st}$ imply that for every $z$ the endomorphism $u$ of $\R^2$
defined by  $u(z):= - (J_{st} + J_\delta(Z))^{-1}(J_{st} - J_\delta(Z))$ is anti
$\C$-linear that is $ u \circ J_{st} = -J_{st} \circ u$. Thus  $u$ is a
composition of the complex conjugation
and a $\C$-linear operator. Denote by $Q_J(z)$ the complex $n \times n$ matrix
such
that $u(z)(v) = Q_J(z)  \overline v$ for any $v \in \C^n$.
 The entries of the matrix $Q_J(z)$
are smooth functions of $z$ and $Q_J(0) = 0$.
The
$J$-holomorphy condition $J(z) \circ dz = dz \circ J_{st}$ can be
written in the form:

\begin{eqnarray}
\label{CR}
z_{\overline\zeta} + Q_J(z) \overline{z}_{ \overline \zeta}   = 0
\end{eqnarray}

Similarly to the proof of lemma \ref{lemma1} consider the isotropic dilations
$d_{\lambda}$. Since the structures $J_\lambda:= (d_\delta)_*(J)$
converge to $J_{st}$ in any $C^\alpha$-norm  as
$\lambda \longrightarrow 0$, we have $Q_{J_\lambda} \longrightarrow 0$ in any
$C^\alpha$ norm. Thus, shrinking $U$ if necessary and using the isotropic
dilations of coordinates as in the proof of lemma \ref{lemma1}, we can assume
that for given $\alpha > 0$ we have $\parallel Q_J
\parallel_{C^\alpha}< <1$ on the unit ball of $\C^n$. In particular, the
system (\ref{CR}) is elliptic.

According to  classsical results \cite{Ve}, the Cauchy-Green
 transform

$$T_{CG}(f) = \frac{1}{2\pi i} \int\int_{\mathbb D} \frac{f(\tau)}{\zeta -
  \tau}d\tau \wedge d\overline\tau$$
is a continuous linear operator from $C^\alpha(\overline{\mathbb D})$ into
 $C^{\alpha+1}(\overline{ \mathbb D})$ for any non-integral $\alpha > 0$ .
Hence  the operator
$$\Psi_{J}: z \longrightarrow w =  z + T_{CG}Q_{J}(z) \overline {z}_{\overline \zeta} $$
takes the space   $C^{\alpha}(\overline{\mathbb D})$  into itself and we can write   the
equation (\ref{CR}) in the form
$[\Psi_J(z)]_{\overline \zeta}  = 0$. Thus, the disc $z$ is $J$-holomorphic if
and only if the map $\Psi_{J}(z):\mathbb D \longrightarrow \C^n$ is
$J_{st}$-holomorphic.
If the norm of $Q_J$ (that is the initial neighborhood $U$) is small enough,
then  by the implicit function theorem the operator    $\Psi_J$
is invertible  and we obtain a bijective
correspondence between {\it small enough} $J$-holomorphic discs and usual
holomorphic discs. This easily implies the existence of a $J$-holomorphic disc
in a given tangent direction through a given point, smooth dependence of such a
disc  on the initial data , as well as the interior elliptic regularity of
discs (the Nijenhuis-Woolf theorem \cite{NiWo}, see
\cite{Si} for further details).

\section{Levi form and plurisubharmonic functions}
Recall some standard definitions playing a substantial role in the sequel .
Let  $r$ be a $C^2$ function on $(M,J)$. We denote by $J^*dr$ the
 differentail form acting on a vector field $X$ by $J^*dr(X):= dr(JX)$.
For example, if $J = J_{st}$ on $\R^2$, then $J^*dr = r_ydx - r_xdy$.
The value of {\it the  Levi
form of} $r$ at a point $p \in M$ and a vector $t \in T_p(M)$ is defined  by
$$L^J_r(p;t):=
-d(J^* dr)(X,JX)$$ where $X$ is an arbitrary smooth vector field in a
neighborhood of $p$ satisfying $X(p) = t$.  This definition is independent of
the choice of vector fields. For instance, if $J = J_{st}$ in $\R^2$, then
$-d(J^*dr) = \Delta r dx \wedge dy$ ($\Delta$ denotes the Lapacian). In
particular, $L_r^{J_{st}}(0,\frac{\partial}{\partial x}) = \Delta r(0)$.

The following properties of the Levi form are fundamental:

\begin{proposition}
\label{pro1}
Let $r$ be a real function of class $C^2$ in a neighborhood of a point $p \in M$.
\begin{itemize}
\item[(i)] If $F: (M,J) \longrightarrow (M', J')$ is a $(J,
  J')$-holomorphic map,  and $\varphi$ is a real function of class
  $C^2$ in a neighborhood of $F(p)$, then for any $t \in T_p(M)$ we have
$L^J_{\varphi \circ F}(p;t) = L^{J'}_\varphi(F(p),dF(p)(t))$.
\item[(ii)] If $z:\mathbb D \longrightarrow M$ is a $J$-holomorphic disc satisfying
  $z(0) = p$, and $dz(0)(e_1) = t$ (here $e_1$ denote the vector
  $\frac{\partial}{\partial \Re \zeta}$ in
  $\R^2$), then $L^J_r(p;t) = \Delta (r \circ z) (0)$.
\end{itemize}
\end{proposition}
 Property (i) expresses the invariance of the Levi form with
respect to biholomorphic maps. Property (ii) is often useful in order to compute the Levi form if a vector
$t$ is given.
\proof (i) Since the map $F$ is $(J, J')$-holomorphic, we have ${J'}^* dr(dF(X))
= dr(J' dF(X)) = dr ( dF(J X)) = d(r \circ F)(JX)$ that is
$F^*({J'}^* dr) = J^* d(r \circ F)$. By the invariance of the exterior
derivative we obtain that $F^*(d{ J'}^* dr) = d J^* d (r \circ F)$. Again
using the holomorphy of $F$, we get $d{ J'}^* dr(dF(X),J'dF(X)) =
F^*(d{ J'}^* dr)(X,JX) = d J^* d(r \circ F)(X,JX)$ which
implies (i).

(ii) Since $z$ is a
$(J_{st},J)$-holomorphic map, (i) implies that $L_r^J(p,t) = L_{r \circ
  z}^{J_{st}}(0,e_1) = \Delta(r \circ z)(0)$. This proves proposition.

\bigskip

As usual, we say that an upper semicontinuous function $r$ on $(M,J)$ is
{\it plurisubharmonic} if its composition with any $J$-holomorphic disc is
subharmonic on $\mathbb D$. For a $C^2$ function this is, in view of proposition
\ref{pro1}, equivalent to the
positive semi-definiteness of the Levi form:
$$L_r^J(p,t) \geq 0 \mbox{ for any } p \in M \mbox{ and } t \in T_p(M) $$
We use
the standard notation $PSH(M)$ for the class of plurisubharmonic functions in $M$.
We say that a
$C^2$ function $r$ is {\it strictly plurisubharmonic} on $M$, if $L_r^J(p,t) > 0$
for any $p \in M$ and $t \in T_p(M) \backslash \{ 0\}$.\smallskip\\

Since the equations (\ref{CR}) are written in complex coordinates, it is
convenient to define the Levi form on the complexified bundle
$T^{(1,0}(M)$. Recall that the bundles $T(M)$ and $T^{1,0}(M)$ are
canonically isomorphic by $X \mapsto V = X - iJX \in T^{1,0}(M)$ for any
section $X$ of $T(M)$ or
equivalently $X = 2 \Re V \in T(M)$ for any section $V$ of $T^{1,0}(M)$.
So we simply define the Levi form of $r$ on $T^{1,0}(M)$ by $L_r^J(p;V) =
-d(J^*dr)(2\Re V, 2J\Re V)(p)$ for any smooth vector field $V$ on
$T^{1,0}(M)$.
\bigskip

Our approach is based on the observation that the Levi form of a function $r$ at
a point $p$ in
an almost complex manifold $(M,J)$ coincides with the Levi form with respect
to the standard structure $J_{st}$ of $\R^{2n}$ if {\it suitable} local
coordinates near $p$ are choosen. Let us explain how to construct these adapted
coordinate systems.

As above, choosing local coordinates near $p$ we may identify a neighborhood
of $p$ with a neighborhood of the origin and assume that $J$-holomorphic discs
are solutions of (\ref{CR}).

\begin{lemma}
\label{normalization}
There exists a  second order polynomial local diffeomorphism $\Phi$ fixing the
origin and with linear part equal to the identity such that in the new coordinates
 the matrix function $Q$ (we drop the index $J$) from the equation (\ref{CR}) satisfies
$$Q(0) = 0, Q_{z}(0) = 0$$
\end{lemma}
Thus, by a suitable local change of coordinates one can remove the linear
terms in $z$ in the matrix $Q$. We stress that in general it is impossible to
get rid of  first order terms containing  $\overline z$ since this would
impose a restriction on the Nijenhuis tensor $J$ at the origin.

The analog of this statement is well-known in the classical theory of
elliptic systems on the complex plane \cite{Ve}. In the
above form the statement  appeared first in E.Chirka's notes \cite{Ch}. In
\cite{SuTu} it is shown that, in an almost complex manifold of (complex) dimension 2,
such a normalization is possible along a given embedded $J$-holomorphic disc;
the proof requires a solution of some $\overline\partial$-type problems. Since the
present assertion provides the normalization only at a given point, the proof is much simpler and works in
any dimension. For convenience of the reader we include the proof following \cite{SuTu}.

\proof Set $z' = \Phi(z)$ and $J' := \Phi_*(z)$.  The $J'$-holomorphy
equations for the disc $z'$ are similar to (\ref{CR}) with the matrix $Q'$
instead of $Q$. We need to establish a relation between the matrices $Q$ from (\ref{CR})
and $Q'$. We have
\begin{eqnarray*}
z'_\zeta = (-\Phi_z Q + \Phi_{\overline\zeta})\overline{z}_{\overline\zeta}, \,\,\,\overline{z'_\zeta} = (\overline{\Phi}_{\overline\zeta} -
\overline{\Phi}_\zeta Q)\overline{z}_{\overline\zeta}
\end{eqnarray*}
Substituting these expressions to the $J'$-holomorphy equation for $z'$, we
obtain the condition $N(z) \overline{z}_{\overline\zeta} = 0$ with
$$N(z) = (-\Phi_z Q + \Phi_{\overline\zeta}) + Q'
(\overline{\Phi}_{\overline\zeta} - \overline{\Phi}_z Q)$$
Since by the Nijenhuis-Woolf theorem for every point $q$ and every vector $v \in
T^{1,0}_q(M)$ there exists a solution $z$ of (\ref{CR}) satisfying $z(0) = q$,
$dz(0)(\frac{\partial}{\partial \zeta}) = v$, we obtain $N = 0$, i.e.
\begin{eqnarray*}
Q' = (\Phi_z Q - \Phi_{\overline z}) (\overline{\Phi}_{\overline z} -
\overline{\Phi}_z Q)^{-1}
\end{eqnarray*}
Set $\Phi(z) = z + \sum_{k,j=1}^{n}\phi_{kj}z_k\overline{z}_j$, where $\phi_{kj}$
are vectors in $\C^n$ with the entries $(\phi_{kj}^s)_{s=1}^n$. Then $\Phi_z =
\sum_{k=1}^n \Phi_k z_k$ and  the matrices $\Phi_k \in M_n(\C)$ have the
entries $(\phi_{kj}^s)_{j,s=1}^n$. Furthermore, $\Phi_z = I + O(\vert z
\vert)$, where $I$ is the unit $n \times n$-matrix.

On the other hand, the Taylor expansion of $Q$ has the form
$Q(z) = A(z) + B(\overline z) + O(\vert z \vert^2)$ where
$A(z) = \sum_{k=1}^n A_k z_k$, $A_k \in M_n(\C)$, and
$B(\overline z) = \sum_{k=1}^n B_k \overline{z}_k$, $B_k \in M_n(\C)$.
Substituting this into the expression for $Q'$, we obtain that
$$Q' = \sum_{k=1}^n (A_k - \Phi_k)z_k + B(\overline z) + O(\vert z \vert^2).$$
Now if we set $\Phi_k = A_k$ (this condition uniquely determines the quadratic
part of $\Phi$), we obtain that $Q'(z'(z))_z(0) = 0$. Since $z' = z + O(\vert
z \vert^2)$, this implies that $Q'_{z'}(0) = 0$ which proves lemma.

\bigskip

The following useful  statement  allows to extend
to the almost complex case all invariance properties of the Levi form without
additional computations. Essentially it follows from the results of  E.Chirka \cite{Ch}.
\begin{proposition}
\label{pro2}
Assume that the local coordinates in a neighborhood of a point $p \in M$ are
chosen according to the previous lemma, that is $Q(0) =0$, $Q_z(0) = 0$ in
(\ref{CR}). Then for any function $r$ of class $C^2$ in a neighborhood of the
origin we have
$$L_r^J(0;t) = L_r^{J_{st}}(0;t)$$
for every $t \in T_0(M)$.
\end{proposition}
\proof According to the Nijenhuis-Woolf theorem, there exists a solution $z$
of (\ref{CR}) such that $z(\zeta) = t\zeta + a\zeta^2 + b\overline{\zeta}^2 +
c\zeta\overline\zeta + o(\vert \zeta \vert^2)$, where $a,b,c \in \C^n$. Since
$Q_z(0) = 0$,we have   $Q(z) = B(\overline z) + o(\vert z \vert)$,  with $B$
as in the proof of the previous lemma. Substituting these expressions into the
equation (\ref{CR}), we obtain
$$2b\overline\zeta + c\zeta +B(\overline t)t\overline\zeta + o(\vert\zeta\vert) = 0$$
so that $c=0$.

Thus, $z(\zeta) = t\zeta + a\zeta^2 + b\overline\zeta^2 + o(\vert \zeta
\vert^2)$. Consider the Taylor expansion of $r$ at the origin:
$$r(z) = r(0) + 2\Re L(z) + 2\Re K(z) + H(z,\overline z) + o(\vert z
\vert^2)$$
where $L$ is a $\C$-linear form, $K$ is a complex quadratic form and $H$ is
the complex Hessian of $r$. Then the terms of degree $\leq 2$ in $\Re L \circ z(\zeta)$ and $\Re K \circ z
(\zeta)$ are harmonic and thus the Laplacian $\Delta (r \circ z)(0)$ coincides with
$L_r^{J_{st}}(0;t)$. On the other hand $\Delta (r \circ z)(0)$ is equal to
$L_r^{J}(0;t)$ by (ii) of proposition \ref{pro1}. This completes the proof.

 \begin{corollary}
For every point $p$ of an almost complex manifold $(M,J)$ there exists a
neighborhood $U$ and local coordinates $z:U \longrightarrow \C^n$ such that
for any function $r$ of class $C^2$ in a neighborhood of $p$ we have $L_r^J(p;
t) = L_{r \circ z^{-1}}^{J_{st}}(0;dz(p)(t))$ for any $t \in T_p(M)$.
\end{corollary}

\section{Bounded strictly plurisubharmonic exhaustion functions}

We begin with   the fundamental definition of Levi convexity.
Let $p$ be a boundary point of a domain $\Omega$ in an almost complex manifold
$(M,J)$; assume that $b\Omega$ is of class $C^2$ in a neighborhood $U$ of $p$.
Then $\Omega \cap U = \{ q \in U: r(q) < 0 \}$ where $r$ is a real function of
class $C^2$ on $U$, $dr(p) \neq 0$.
\begin{definition}
$\Omega$ is called Levi convex  at $p \in b\Omega$
if $L_r^J(p;t) \geq 0$ for any $t \in T_p(b\Omega) \cap J(T_p(b\Omega))$ and  strictly Levi convex
 at $p$  if  $L_r^J(p;t) > 0$ for any non-zero $t \in
T_p(b\Omega) \cap J(T_p(b\Omega))$.  If  $\Omega$ is a relatively compact
domain with $C^2$ boundary in an almost complex manifold $(M,J)$, then $\Omega$
is called Levi convex  if it is Levi convex at every boundary
point.
\end{definition}
It is easy to show that this
definition does not depend on the choice of defining functions. In the works
\cite{GrEl,ElTh,Mc} strictly Levi convex domains are called $J$-convex. We
prefer the terminology closer to traditional complex analysis.

Recall also that a  continuous  map  $u: \Omega \longrightarrow [a,0[
\subset \R$ is called a {\it bounded exhaustion function} for $\Omega$ if for every
$a \leq b < 0$ the pull-back $u^{-1}([a,b])$ is compact in $\Omega$.

Our first main result is the following
\begin{theorem}
\label{theo1}
Let $(M,J)$ be an almost complex manifold and let $\Omega
\subset M$ be a relatively compact Levi convex domain  with $C^3$
boundary, such that there exists a $C^2$ strictly plurisubharmonic
function $\phi$  in a neighborhood $U$ of $b\Omega$. Let $r$ be any
$C^3$ defining function for $\Omega \cap U$. Then there exists a neighborhood $U'$ of
$b\Omega$ and constants $A > 0$, $0 < \eta_0 < 1$, such that for any $0 <
\eta \leq \eta_0$ the function $\rho = -(-re^{-A\phi})^\eta$ is strictly
plurisubharmonic on $\Omega \cap U'$. If $U$ is a neighborhood of
$\overline \Omega$, then $\rho$ is strictly plurisubharmonic on $\Omega$.
\end{theorem}

Thus, $\rho$ is a bounded strictly plurisubharmonic exhaustion function for
$\Omega$. The proof is based on the method of K.Diederich - J.E.Fornaess \cite{DiFo1} (slightly
modified
by M.Range in \cite{Ra}). In its first step the following expression for the
Levi form of $\rho$ is determined.

\begin{lemma}
\label{Leviform}
Under the hypothesis of  Theorem \ref{theo1} there exists a neighborhood $U'$ of
$b\Omega$ such that for every $p \in \Omega \cap U'$ and
$v \in T_p(M)$ we have
\begin{eqnarray*}
L_\rho^{J}(p;v) = \eta(-r)^{\eta - 2}e^{-\eta A \psi}D(v)
\end{eqnarray*}
where
\begin{eqnarray*}
& &D(v) = Ar^2(p)[L_\psi^{J}(p;v) - \eta A \vert
\partial_{J}\psi(p)(v)\vert^2] + (- r(p))[L_r^{J}(p;v) - 2 \eta A \Re
\partial_J r(p)(v)\overline{\partial_{J}\psi(p)(v)}] \\
& &+ (1 - \eta)\vert \partial_{J}r(p)(v)\vert^2
\end{eqnarray*}
\end{lemma}
In the case of the standard structure of $\C^n$ this formula is due to \cite{DiFo1}.
\proof Fix a point $p \in \Omega$ close enough to $b\Omega$ and choose
coordinates $z:W \longrightarrow \C^n$ in a neighborhood $W$ of $p$  according to lemma
\ref{normalization}. In these coordinates
the neighborhood $V$ can be identified with a neighborhood $z(V)$ of the
origin in $\C^n$. For simplicity of notation, we denote again by $J$ the
direct image $z_*(J)$, and by $\rho$ (resp. $r$) the composition $\rho \circ
z^{-1}$ (resp. $r \circ z^{-1}$). In particular, we have $J(0) = J_{st}$. Consider a vector
 $t \in T_0^{1,0}(\C^n)$. According  to \cite{DiFo1} we have the following
 expression
 for the Levi form of the function $\rho = -(-r e^{-A\psi})^\eta$ with respect
 to the standard structure $J_{st}$:

\begin{eqnarray*}
L_\rho^{J_{st}}(0;t) = \eta(-r)^{\eta - 2}e^{-\eta A \psi}D_{st}(t)
\end{eqnarray*}
where
\begin{eqnarray*}
& &D_{st}(t) = Ar^2(p)[L_\psi^{J_{st}}(0;t) - \eta A \vert
\partial_{J_{st}}\psi(0)(t)\vert^2] + ( - r(p))[L_r^{J_{st}}(0;t) - 2 \eta A \Re
\partial_{J_{st}} r(0)(t)\overline{\partial_{J_{st}}\psi(0)(t)}] \\
& &+ (1 - \eta)\vert \partial_{J_{st}}r(0)(t)\vert^2
\end{eqnarray*}

On the other hand, by proposition \ref{pro2} for every real function $\varphi$
of class $C^2$ in a neighborhood of the origin its   Levi form at $0$ with respect to $J$
coincides with the Levi form with respect to $J_{st}$. Furthermore, the condition $J(0) = J_{st}$ implies
$\partial_{J_{st}}\varphi(0) = \partial_J \varphi(0)$. Therefore, for any
vector $v \in T_p^{1,0}(M)$ we have a similar expression for the Levi form
$L_\rho^J(p;v)$ with respect to $J$ without any additional calculations. This
proves lemma.

\bigskip

Now the proof of Theorem \ref{theo1} is quite similar to \cite{DiFo1,Ra}. For
the sake of completeness we include the details.

{\it Proof of Theorem \ref{theo1}:}
Recall that any
almost complex manifold admits a Hermitian metric (see, for instance,
\cite{Ko}, Vol. II). Fix such a metric $d\mu^2$ on $M$ and denote by $\parallel
. \parallel_p$ the norm  induced by $d\mu^2$ on $T^{1,0}_p(M)$. Then we have the
orthogonal decomposition $T_p^{1,0}(M) = T^t_p \oplus T^n_p$ where $T^t_p = \{
v \in T_p^{1,0}(M): \partial r_J(p)(v) = 0 \}$ is the ``tangent'' space and its
orthogonal complement $T^n_p$ is the ``normal'' space. If $p \in b\Omega$,
then $T^t_p$ is canonically isomorphic to the holomorphic tangent space
$T_p(b\Omega) \cap J(T_p(b\Omega))$  by the canonical identification of
$T_p(M)$ and $T_p^{1,0}(M)$. So for any vector $v
\in T_p^{1,0}(M)$ we have the decomposition $v = v^t + v^n \in T^t_p \oplus
T^n_p$ into the ''tangent'' and ``normal'' components. Fix a neighborhood $U$ of
a point $q \in b\Omega$ of the form $U = b\Omega \times ]-\varepsilon,\varepsilon[$ and
denote by $\pi: W \longrightarrow b\Omega$ the natural projection.
Since $r$ is of class
$C^3$, given smooth vector field $V: p \longrightarrow V(p) \in T_p^{1,0}(M)$ the function
$p \mapsto L_r^J(p,V^t(p))$ is of class $C^1$  on $W$. Furthemore, the
function $p \mapsto L_r^J(p;V^t(p)) - L_r^J(\pi(p); (d\pi(p)(V(p)))^t)$ vanishes on
$b\Omega$. In local coordinates the last expression is just a quadratic form
in $V$ with coefficients of class $C^1$ in $p$ vanishing on $b\Omega$. Therefore, there exist a constant $C_1 > 0$ such that
$$\vert L_r^J(p;V^t(p)) - L_r^J(\pi(p);(d\pi(p)(V(p)))^t) \vert \leq C_1 \vert r(p) \vert
\parallel V(p) \parallel_p^2.$$
Since $\Omega$ is Levi convex, we have $L_r^J(p;(d\pi(p)(V(p)))^t) \geq 0$ for every
$V$. This implies
\begin{eqnarray}
\label{1}
L_r^J(p;V^t(p)) \geq -C_1\vert r(p) \vert \parallel V(p) \parallel_p^2
\end{eqnarray}
By compactness argument we can assume that there exists a neighborhood $W$ of
$b\Omega$ such that the estimate (\ref{1}) holds
for every $p \in \Omega \cap W$  and any $(1,0)$ vector field
$V$ on $W$.

If $V_1$ and $V_2$ are $(1,0)$ vector fields, then we set
$L_r^J(p;V_1,V_2) = -d(J^*dr)(2\Re V_1, 2J\Re V_2)$. We have
\begin{eqnarray*}
& &\vert L_r^J(p;V(p)) - L_r^J(p;V^t(p))\vert = \vert L_r^J(p;V^t(p),V^n(p)) +
L_r^J(p;V^n(p),V^t(p))\\
& & + L_r^J(p;V^n(p),V^n(p))\vert \leq C_2 \parallel V(p)
\parallel_p \parallel V^n(p) \parallel_p
\end{eqnarray*}
and
\begin{eqnarray*}
\parallel V(p) \parallel_p \leq C_3 \vert \partial_J r(V(p))\vert
\end{eqnarray*}
Together with (\ref{1}) these estimates imply that
\begin{eqnarray}
\label{2}
L_r^J(p;V(p)) \geq -C_4 \vert r(p) \vert \parallel V(p) \parallel_p^2 - C_4
\parallel V(p) \parallel_p \vert \partial_J r(V(p)) \vert
\end{eqnarray}
for any $ p \in \Omega \cap W$ and any vector field $V$.

Now we can estimate from above $D(v)$.
Fix $C_5 > 0$ such that $L_\psi^J(q;v) \geq C_5\parallel v \parallel_p^2$ for
any $q \in W \cap \Omega$ and
any $v \in T_q^{1,0}(M)$. Then there exists $\eta_0(A)$ such that for any $0 < \eta <
\eta_0(A)$

\begin{eqnarray*}
D(v) \geq Ar^2(C_5 - C_5/2)\parallel v \parallel^2 + (- r) (r C_4 \parallel v
\parallel^2 - C_6 \vert \partial_J r(p)(v)\vert \parallel v \parallel) +
(1/2)\vert \partial_Jr(p)(v) \vert^2
\end{eqnarray*}
Since
\begin{eqnarray*}
C_6 \vert r(p) \vert  \vert \partial_J r(p)(v) \vert \parallel v
\parallel \leq (1/4)\vert \partial_J r(p)(v)  \vert^2 + C_7 r^2
\parallel v \parallel^2
\end{eqnarray*}
we obtain
\begin{eqnarray*}
D(v) \geq r^2(A C_5/2 - C_8)\parallel v \parallel^2
\end{eqnarray*}
where all constants $C'$s are positive and independent of $A$, $\eta$ and
$v$. Now this is enough to fix $A > 2C_8/C_5$ and then $\eta_0 = \eta(A)$.

Finally, let $\psi$ be strictly plurisubharmonic on $\overline \Omega$. Then
$L_\psi^J(p;v) \geq C_9 \parallel v \parallel^2$ for any $p \in
\overline\Omega$ and $v \in T^{1,0}_p(M)$ with
$C_9 > 0$; moreover, $K: = \Omega \backslash W$ is a compact subset in
$\Omega$
and there exists a $\delta > 0$ such that $r(p) \geq \delta$ for every $p \in
K$. The expression for $D(v)$ given by lemma \ref{Leviform} now holds in
$\Omega$.  Then there for every $0 < \eta < \eta_1(A)$  we have
 $D(v) \geq (A\delta^2C_9/2  - C_{10})\parallel v
\parallel^2$ for all $p \in K$ and $v \in T^{1,0}_p(M)$.  So we take $A >
\max(2C_8/C_5,2C_{10}/C_9\delta^2)$ and then fix $\eta$. This proves theorem. \smallskip \\

\begin{remark}
1) As in the case of $\C^n$, for a fixed point $p \in
b\Omega$ and any given $0 < \eta < 1$ there is a neighborhood $U$ of $p$, a strictly
plurisubharmonic function $\psi$ in $U$ and $A > 0$ such that $\rho$ is
strictly plurisubharmonic in $\Omega \cap U$. Indeed, in local coordinates
given by lemma \ref{normalization} consider the function $\psi(z) =
\sum_{j=1}^n \vert z_j \vert^2$. Then $\psi$ is strictly plurisubharmonic in $U$ and $d\psi(p) = 0$
so  the result follows from the previous estimate of $D(v)$. \smallskip \\
2) We want to remind the reader of the fact, known already for the classical situation,
that, in general, a bounded strictly plurisubharmonic exhaustion for $\Omega$ cannot
be chosen to be more than just Hölder continuous up to $b\Omega$ (see \cite{DiFo1}.
In fact, the infimum of all possible $0<\eta <1$ in Theorem \ref{theo1} is an interesting
invariant of $\Omega$ linked to the $\overline\partial$-Neumann problem in
classical complex analysis (see \cite{Koh}).
\end{remark}

As an obvious consequence we obtain the following

\begin{corollary}
No $J$-holomorphic discs can touch $b\Omega$ from the inside.
\end{corollary}
In the case where $\Omega$ is strictly Levi convex, this statement is
well-known (see, for instance, \cite{Mc}, lemma 2.4).

\section{Characterization of Stein structures: The continuity principle}

Let $\Omega$ be a relatively compact domain in an almost complex manifold $(M,J)$. Let $D:
\mathbb D \longrightarrow M$ be a $J$-holomorphic disc continuous on
$\overline{\mathbb D}$. With some abuse of notation we will denote its image
$D(\mathbb D)$ just by $D$, $\overline D = D(\overline{\mathbb D})$  and denote by $bD$ its ``boundary'': $bD: =
D(b\mathbb D)$.  By a {\it
Hartogs family }  we mean  a continuous map $D: \overline{\mathbb D} \times [0,1]
\longrightarrow M$ such that for every $t \in [0,1]$ the map $D_t:= D(\bullet,t)$ is
$J$-holomorphic on $\mathbb D$, for every $t \in ]0,1]$ we have $\overline D_t \subset
\Omega$ and $bD_0 \subset \Omega$.
\begin{definition}
We say that  $\Omega$ is
  disc-convex  if for any Hartogs family of discs we have $D_0 \subset
\Omega$.
\end{definition}
The classical definition of pseudoconvexity also can be
extended to the almost complex case without changes.
\begin{definition}
A domain $\Omega$ in an almost complex manifold $(M,J)$ is pseudoconvex if for
any compact subset $K$ in $\Omega$ its plurisubharmonically convex hull
$\hat K_\Omega := \{ p \in \Omega: u(p) \leq \sup_{q \in K} u(q) , u \in PSH(\Omega)
\cap C(\Omega) \}$ is compact.
\end{definition}

Finally, recall that a continuous proper map  $u:\Omega \longrightarrow [0,+\infty[$ is
called an {\it exhaustion function} for a domain $\Omega$.

\begin{definition}
A domain $\Omega$ in an almost complex manifold $(M,J)$ is a Stein domain if
there exists a strictly plurisubharmonic exhaustion function on $\Omega$.
\end{definition}

The following consequence of Theorem \ref{theo1} is the almost complex analog of the results
valid in the classical situation.

\begin{theorem}
\label{theo2}
Let $\Omega$ be a relatively compact domain with boundary of class $C^3$
in an almost complex manifold $(M,J)$. Suppose that $M$ admits a strictly
plurisubharmonic function. Then the following conditions are equivalent:
\begin{itemize}
\item[(i)] $\Omega$ is Levi convex;
\item[(ii)] $\Omega$ is disc-convex;
\item[(iii)] $\Omega$ admits a bounded strictly plurisubharmonic exhaustion
  function;
\item[(iv)] $\Omega$ is pseudoconvex;
\item[(v)] $\Omega$ is a Stein domain.
\end{itemize}
\end{theorem}

\proof $(i) \Longrightarrow (iii)$ by Theorem \ref{theo1}. The proof of $(iii)
\Longrightarrow (ii)$ is quite
similar to the case of $\C^n$. Next, we prove that $(ii) \Longrightarrow (i)$. Suppose that (ii) holds and $\Omega$ is
not Levi convex. Then there exists a boundary point $p \in b\Omega$ and a
vector $v \in H_p^J(b\Omega)$ such that $L_r^J(p,v) < 0$, where $r$ is a local
defining function of $\Omega$ near $p$. This condition is open and stable with respect
to small translations of $b\Omega$ along the normal direction at $p$, so the
same holds for the hypersurfaces $b\Omega_t = \{ r = -t \}$ if $t \geq 0$ is small
enough. It follows from results of J.-F.Barraud - E.Mazzilli  \cite{BaMa}
and S.Ivashkovich - J.P.Rosay \cite{IvRo} that for every $t$ there exists a $J$-holomorphic disc
$D_t$ such that $(\rho \circ D_t)(0) = -t$ and $(\rho \circ D_t)(\zeta) < -t$
if $\zeta \in \overline{\mathbb D} \backslash \{ 0 \}$. This family depends
continuously on $t$ and so this is a Hartogs family such that $D_0$ is not
contained in $\Omega$: a contradiction.  So (i), (ii) and (iii) are
equivalent. The  proofs of  $(iii) \Longrightarrow (iv)
\Longrightarrow (ii)$
and $(iii) \Longrightarrow (v) \Longrightarrow (iv)$ are
quite similar to the classical arguments  for the standard $\C^n$.  This
completes the proof.

\bigskip

Y.Eliashberg and M.Gromov \cite{GrEl}  proved
that every Stein domain admits the canonical symplectic structure defined
by the symplectic form $\omega_u = -dJ^*du$, where $u$ is a strictly
plurisubharmonic exhaustion function. Clearly, $\omega_u(v,Jv) > 0$ for any
non-zero tangent vector $v$ that is $J$ is tamed by the symplectic form
$\omega$ in the sense of Gromov. Furthermore, if $\tilde u$ is another strictly
plurisubharmonic exhaustion function for $\Omega$, then the symplectic
manifolds $(\Omega,\omega_u)$ and $(\Omega,\omega_{\tilde u})$ are
symplectomorphic \cite{GrEl}. In this sense the above symplectic structure
defined by a strictly plurisubharmonic  exhaustion function is
canonical.  Thus, we have the following

\begin{corollary}
Let $(M,J)$ be an almost complex manifold admitting a strictly
plurisubharmonic function. Then a relatively compact  domain $\Omega$
with $C^3$ boundary in $M$ admits the canonical symplectic structure if and only
if $\Omega$ is pseudoconvex.
\end{corollary}

\section{Normal families and taut manifolds}
Let $(M,J)$ be an almost complex
manifold. Recall that it is called {\it taut} if any sequence $(D_k)$ of
$J$-holomorphic discs in $M$  either contains a compactly convergent
subsequence or is compactly
divergent. This is well-known and easily follows from the standard elliptic
estimates for the Cauchy-Green kernel (see for instance, \cite{Si}) showing that the
limit in the compact open topology of a sequence of $J$-holomorphic discs also
is a $J$-holomorphic disc. Since Gromov's compactness theorems for $J$-holomorphic discs can be
viewed as normal family type theorems, the class of taut  almost complex manifolds
is very natural. It is well-known that  any complete hyperbolic manifold is
taut; the inverse in general is not true.

\begin{theorem}
\label{theo3}
Suppose that an almost complex manifold  $(M,J)$   admits a strictly
plurisubharmonic function. Then any relatively compact pseudoconvex domain with $C^3$
boundary in $M$  is taut.
\end{theorem}
In the case of $\C^n$ this is an immediate corollary of the K.Diederich -
J.E. Fornaess
result \cite{DiFo1,DiFo2} on the existence of a bounded strictly plurisubharmonic exhaustion
function
and the Montel theorem. In the case of an arbitrary complex manifold (with an
integrable structure) the result is due to N.Sibony \cite{Sib}. We also point
out that if in the hypothesis of Theorem \ref{theo3} the boundary of the
domain is strictly Levi convex,
then the domain is complete hyperbolic \cite{GaSu,IvRo}.

\proof Since $M$ admits a strictly plurisubharmonic function, it follows from
the results of \cite{GaSu, IvRo}  that $M$ is hyperbolic at every point (in
the sense of H.Royden \cite{Ro}) and so
hyperbolic for the Kobayashi distance by the classical H.Royden theorem \cite{Ro} (see, for instance,
\cite{IvRo,Kr} for the almost complex version of this theorem). In particular,
a domain $\Omega$ satisfying the hypothesis of Theorem \ref{theo3} is a
hyperbolic manifold and the Kobayashi distance determines the usual topology
on $\Omega$ (see, for example, \cite{Ko}).
By Theorem \ref{theo1} $\Omega$ admits a bounded strictly
plurisubharmonic exhaustion function $\rho$. Now given a family of
$J$-holomorphic discs $(D_k)$ in $\Omega$ it suffices to apply the
argument of  the proof of corollary 5 in \cite{Sib} since it only uses the
subharmonicity of $\rho \circ D_k$  and so does not require any modifications.

\section{Approximation of confoliations by contact structures}

According to Y.Eliashberg - W.Thurston \cite{ElTh} a tangent hyperplane field
$\xi = \{ \alpha = 0 \}$  on a $(2n+1)$-dimensional manifold $\Gamma$ is
called a {\it positive confoliation} if there exists an almost complex
structure $J$ on the bundle $\xi$ such that
\begin{eqnarray*}
d\alpha(X,JX) \geq 0
\end{eqnarray*}
for any vector $X \in \xi$. The 1-form $\alpha$ is defined up to the
multiplication by a nonvanishing function. Thus, the confoliation condition for $\xi$ is
equivalent to the existence of a compatible Levi convex CR-structure
(in general, non-integrable). In other words, if $\Gamma = \{ r = 0 \}$ is a
smooth Levi convex  hypersurface in
an almost complex manifold $(M,J)$, then the distribution of its holomorphic
tangent spaces $\xi = T\Gamma \cap J(T\Gamma)$ is a confoliation: we can set
$\alpha = J^*dr$. In particular, if $\Gamma$ is a strictly Levi convex
hypersurface, then $\xi = \{ J^*dr = 0 \}$ is a contact structure. Recall that
a tangent hyperplane field $\xi = \{\alpha = 0 \}$ on $\Gamma$ is called a
contact structure if $\alpha \wedge (d\alpha)^n \neq 0$ on $\Gamma$. One of the main
questions
considered by Y.Eliashberg - W.Thurston concerns the possibility to  deform a
given confoliation to
a contact structure or approximate it by contact structures.
Combining the contact topology techniques with the geometric  foliation theory
 they
obtained several results of this type in the case where $\Gamma$ is of real
dimension 3. Our next result works in any dimension.

\begin{theorem}
\label{theo4}
Let $\Omega$ be a relatively compact pseudoconvex domain with $C^\infty$
boundary in an almost complex manifold $(M,J)$. Assume that there exists a
$C^\infty$ strictly plurisubharmonic function $\psi$ in a neighborhood of $b\Omega$.
Then the confoliation of holomorphic tangent spaces $T(b\Omega) \cap
J(T(b\Omega))$ can be approximated in any $C^{k}$ norm by contact
structures.
\end{theorem}

\proof Fix a  small enough neighborhood $U$ of $b\Omega$ such that $r$ is a
defining function of $\Omega \cap U$ and for some $0
< \eta < 1$ and $A > 0$ the function $\rho = -(-r e^{-A\psi})^{\eta}$ is strictly
plurisubharmonic on $\Omega \cap U$ by Theorem \ref{theo1}. Set also $\varphi = r e^{-A\psi}$. For
 $\delta > 0$ small enough consider the hypersurfaces $\Gamma_\delta = \{ \rho
= - \delta \} = \{ \varphi = - \delta^{1/\eta} \}$. By the Sard theorem the
set of critical values of $\rho$ is of the Lebesgue measure zero, so we can
fix a sequence $(\delta_j)_j$ of non-critical values of $\rho$ converging to
$0$. Then every $\Gamma_j: = \Gamma_{\delta_j}$ is a strictly Levi convex
hypersurface and admits the canonical contact structure $\xi_j$ defined by the form
$J^*d\rho$.  Since a contact form is defined up to multiplication by a
non-vanishing function and $d\rho = \eta (-\varphi)^{\eta - 1}d\varphi$, we
have $\xi_j = \{ J^*d\varphi = 0 \}$. Fix a smooth Riemannian metric on
$M$. Shrinking the neighborhood $U$ if necessary we can assume that it is
foliated by the normals to $b\Omega$. For any $j$ consider a smooth
diffeomorphism $\pi_j:b\Omega \longrightarrow \Gamma_j$ which associates to
every point of $p \in b\Omega$ the point of $\Gamma_j$ defined by the intersection of
$\Gamma_j$ with the normal  of $b\Omega$ at $p$.  Then the pull-back
$\alpha_j: = \pi_j^*(J^*d\varphi)$ defines a contact structure $\xi^*_j =
\{ \alpha_j = 0 \}$ on $b\Omega$ and the sequence of forms $(\alpha_j)_j$
converges to the form $J^*d\varphi$ on $b\Omega$ in any $C^k$ norm which
proves theorem.

\begin{remark}
In a suitable neighborhood of a fixed point $p \in b\Omega$
there always exists a smooth strictly plurisubharmonic function. So
every  confoliation can be approximated locally by contact structures.
\end{remark}

\end{document}